\documentclass[a4paper,12pt]{article}
\usepackage[english]{babel}
\usepackage{amsfonts,amssymb,amsmath}

\begin{document}

 \centerline{\bf\Large  Asymptotics of the convolution of the Airy function} 
\vskip 1mm 
 \centerline{\bf\Large  and a function of the power-like behavior}

 \vskip 5mm 
 \centerline{\bf S.V.~Zakharov} 

 \vskip 5mm 

\begin{center}
Institute of Mathematics and Mechanics,\\
Ural Branch of the Russian Academy of Sciences,\\
16, S.Kovalevskaja street, 620990, Ekaterinburg, Russia
\end{center}
 
\vskip 5mm 

\textbf{Abstract.}
The asymptotic behavior of the convolution-integral
of a special form of the Airy function
and a function of the power-like behavior
at infinity is obtained. 
The integral under consideration is the solution
of the Cauchy problem for an evolutionary
third-order partial differential equation
used in the theory of wave propagation 
in physical media with dispersion.
The obtained result can be applied to studying
asymptotics of solutions of the KdV equation
by the matching method.


 
\vskip 5mm 
 
\section{Introduction}\label{s1} 
In the present paper,
we study the behavior as $|x|+t\to\infty$
of the convolution 
\begin{equation}\label{us} 
 u(x,t) = \frac{1}{\root 3 \of {3t} } 
 \int\limits_{-\infty}^{+\infty} 
 f(y)\, \mathrm{Ai} \left( \frac{x-y}{\root 3 \of {3t} } \right) \, dy 
 \end{equation} 
of the Airy function 
\begin{equation}\label{af} 
 \mathrm{Ai}\,(x) = 
 \frac{1}{\pi} \int\limits_{0}^{+\infty} 
 \cos \left( \frac{\theta^3}{3}+ x\theta \right) \,
d\theta 
 \end{equation}            
 and a locally Lebesgue integrable function~$f$, 
 which satisfies the following asymptotic 
(in the sense of Poincar\'{e}) relations: 
 \begin{equation}\label{le} 
 f(x) = \sum\limits_{n=0}^{\infty} 
 f^{\pm}_n x^{-n}, 
 \qquad x\to \pm\infty.
 \end{equation} 

 Integral~(\ref{us}), 
 which is understood in this case in the sense
of the limit
 $\lim\limits_{R\to +\infty}\int\limits_{-\infty}^{R}$ 
of a standard Lebesgue integral, is of interest 
as the solution of the Cauchy problem for a third-order
equation with the initial function~$f$: 
 $$ 
 \frac{\partial u}{\partial t} + 
 \frac{\partial^3 u}{\partial x^3} =0, \quad u(x,0) =
f(x), 
 \qquad t\geqslant 0, \quad x\in\mathbb{R}. 
 $$ 
 This equation is used 
in the theory of wave propagation 
in physical media with dispersion~\cite{Gbw} and it
is sometimes called the linearized Korteweg--de Vries (KdV) equation. 
 Notice that the Airy function and
the  Airy transform (convolution of a general form) 
have applications in wave optics, quantum mechanics, 
and in modern laser technologies~\cite{vs,at}. 
 In addition, the investigation of the behavior 
of integral~(\ref{us}) 
is of independent interest for asymptotic analysis,
for example, the obtained result can be applied to studying
asymptotics of solutions of the nonlinear KdV equation
by the matching method~\cite{ib}
in the same way as the large-time asymptotics of solutions
of the heat equation~\cite{hd} is applied 
to the Cauchy problem for a nonlinear parabolic equation~\cite{2ps,zz}.

An asymptotics  can be easily found
 for large values of~$t$ and bounded values of~$x$, 
 if $f(x)$ exponentially 
tends to zero as $x \to \infty$. 
In this  case, it suffices to expand the Airy function
 into  a series in $t^{-1/3}$. 
 However, for more general conditions~(\ref{le}) 
the problem becomes much more difficult.

 \bigskip 
\setcounter{equation}{0}
 
\section{Investigation integral}\label{s2} 
 Following the approach used in paper~\cite{hd}, 
we represent function~(\ref{us}) in the form 
\begin{equation}\label{ud} 
 u(x,t)=U^-_1(x,t) + U^-_0(x,t) +U^+_0(x,t) + U^+_1(x,t), 
 \end{equation} 
 where $$ 
 U^-_1(x,t) = \int\limits_{-\infty}^{-\sigma} \dots, 
 \quad 
 U^-_0(x,t) = \int\limits_{-\sigma}^{0}\dots, 
 \quad 
 U^+_0(x,t) = \int\limits_{0}^{\sigma}\dots, 
 \quad 
 U^+_1(x,t) = \int\limits_{\sigma}^{+\infty}\dots, 
 $$ 
 $$ 
 \sigma = (|x|^3 + t)^{p/3}, 
 \qquad 
 0 < p < 1, 
 $$ 
 and dots denote the integrand from formula~(\ref{us}) 
 together with the factor $(3t)^{-1/3}dy$. 
 
 In the integral $U^-(x,t)$,
 we make the change $y = \theta {\root 3 \of {3t}}$, 
 setting $$ 
 \mu = \frac{\sigma}{{\root 3 \of {3t}}}, 
 \qquad 
 \eta = \frac{x}{{\root 3 \of {3t}}}. 
 $$ 
 Using condition~(\ref{le}) as $x\to -\infty$, 
we obtain 
$$ 
 U^-_1(x,t) = \int\limits_{-\infty}^{-\mu} 
 f(\theta {\root 3 \of {3t}}) \mathrm{Ai}(\eta-\theta)
d\theta = 
 $$ 
 \begin{equation}\label{up1i} 
 = \sum\limits_{n=0}^{N-1} 
 {f^-_n} (3t)^{-n/3} 
 \int\limits_{-\infty}^{-\mu} 
 \theta ^{-n} \mathrm{Ai}(\eta-\theta)
d\theta + 
 \int\limits_{-\infty}^{-\mu} 
 R_N(\theta {\root 3 \of {3t}}) \mathrm{Ai}(\eta-\theta) d\theta, 
 \end{equation} 
 where $ R_N(s) = O(|s|^{-N}).$ 
 From the last relation, we have an estimate 
of the remainder: 
 \begin{equation}\label{ern} 
 \int\limits_{-\infty}^{-\mu} 
 R_N(\theta {\root 3 \of {3t}}) \mathrm{Ai}(\eta-\theta) d\theta 
 = O(\sigma^{-N+1}). 
 \end{equation} 
 
Let us find the dependence of the
integral
 $$ 
 \int\limits_{-\infty}^{-\mu} 
 \theta ^{-n} \mathrm{Ai}(\eta-\theta)
d\theta
 $$ 
on value of $\mu$.
 For $n=0$ and $t\geqslant |x|^{\alpha}$, 
 $2 + p < \alpha < 3$, we have $$ 
 \int\limits_{-\infty}^{-\mu} 
 \mathrm{Ai}(\eta-\theta) d\theta = 
 \int\limits_{0}^{+\infty} 
 \mathrm{Ai}(\eta+\theta' ) d\theta' - \int\limits_{0}^{\mu} 
 \mathrm{Ai}(\eta+\theta' ) d\theta' = 
 $$ 
 \begin{equation}\label{ei0} 
 =\int\limits_{\eta}^{+\infty}\mathrm{Ai}(\eta')
d\eta' 
 -\sum\limits_{r=1}^{N-1} \frac{\mu^r}{r!} \mathrm{Ai}^{(r-1)}(\eta) 
 + O(\sigma^{-\gamma N}), 
 \qquad 
 \sigma\to \infty, 
 \end{equation} 
 where $\gamma = \frac{\alpha}{3p} - 1 > 0.$ 
To obtain (\ref{ei0}), we used Taylor's
formula for the expansion of $\mathrm{Ai}(\eta+\theta' )$ 
in  $\theta'$  and the estimates 
\begin{equation}\label{ms} 
 \mu = O\left(t^{\frac{3p-\alpha}{3\alpha}} \right), 
 \qquad 
 \sigma = O\left(t^{p/\alpha} \right), 
 \qquad 
 \mu = O\left( \sigma^{-\gamma} \right). 
 \end{equation} 
as $\sigma\to \infty$ on the set
$ 
 T_{\alpha} = \{ (x,t)\, : \, x\in\mathbb{R},\
t \geqslant |x|^{\alpha} \}.
$ 
For $n\geqslant 1$, using the regularization method,
 i.e., subtraction of singularities
as $\theta\to -0$, 
 we obtain $$ 
 \int\limits_{-\infty}^{-\mu} 
 \theta ^{-n} \mathrm{Ai}(\eta-\theta)
d\theta = 
 \int\limits_{-\infty}^{-1} \theta ^{-n} 
 \mathrm{Ai}(\eta-\theta) d\theta + 
 $$ 
 $$ 
 + \int\limits_{-1}^{-\mu} \Psi_n(\theta,\eta)
d\theta + 
 \sum\limits_{r=0}^{n-1} \frac{(-1)^{r-n}}{r!} \mathrm{Ai}^{(r)}(\eta) 
 \int\limits_{-1}^{-\mu} \theta ^{r-n}
d\theta, 
 $$ 
 where 
\begin{equation}\label{psin} 
 \Psi_n(\theta,\eta) = \theta ^{-n} \left[ 
 \mathrm{Ai}(\eta-\theta) - \sum\limits_{r=0}^{n-1} \frac{(-\theta)^{r}}{r!}\mathrm{Ai}^{(r)}(\eta) 
 \right], 
 \end{equation} 
 and the sum in the square brackets is 
a partial sum of the Taylor  series 
for $\mathrm{Ai}(\eta-\theta)$ in variable $\theta$. 
 Thus, $$ 
 \int\limits_{-\infty}^{-\mu} \theta ^{-n} 
 \mathrm{Ai}(\eta-\theta) d\theta = 
 \int\limits_{-\infty}^{-1} \theta ^{-n} 
 \mathrm{Ai}(\eta-\theta) d\theta 
 -\frac{\mathrm{Ai}^{(n-1)}(\eta)}{(n-1)!} \ln\mu + 
 $$ 
 \begin{equation}\label{ein} 
 -\sum\limits_{r=0}^{n-2} \mu^{r-n+1} 
 \frac{\mathrm{Ai}^{(r)}(\eta)}{r!(r-n+1)} 
 - \int\limits_{-\mu}^{0}\Psi_n(\theta,\eta)
d\theta 
 +B^-_n(\eta), 
 \end{equation} 
 where 
\begin{equation}\label{bnm} 
 B^-_n(\eta)= \sum\limits_{r=0}^{n-2} 
 \frac{\mathrm{Ai}^{(r)}(\eta)}{r!(r-n+1)} 
 + \int\limits_{-1}^{0}\Psi_n(\theta,\eta)
d\theta. 
 \end{equation} 
 From formulas~(\ref{psin}) we conclude
that $\Psi_n(\theta,\eta)$ has no singularities
as $\theta\to -0$; 
 whence we get
 $$ 
 \int\limits_{-\mu}^{0}\Psi_n(\theta,\eta)
d\theta 
 = \sum\limits_{r=1}^{N-1}\mu^r q_{r}\mathrm{Ai}^{(r)}(\eta) + 
 O(\sigma^{-\gamma N}), 
 $$ 
\begin{equation}\label{pinm} 
 \int\limits_{-1}^{0} 
 \Psi_n(\theta,\eta) d\theta = 
 O\left( |\eta|^{-1/4+n/2} \right), 
 \qquad \eta\to \infty. 
 \end{equation} 
 
 Substituting relations~(\ref{ei0}) and~(\ref{ein}) 
 in formula~(\ref{up1i}), and using estimate~(\ref{ern}), 
 we find 
\begin{equation}\label{up1f} 
 U^-_1(x,t) = 
 f^-_0\int\limits_{\eta}^{\infty} 
 \! \mathrm{Ai}(\theta)\, d\theta 
 +\sum\limits_{n=1}^{N-1} 
 f^-_n (3t)^{-n/3} 
 \left[ \int\limits_{-\infty}^{-1} 
 \theta^{-n} \mathrm{Ai}(\eta-\theta)
d\theta +B^-_n(\eta) 
 \right]+ 
 \end{equation} 
 $$ 
 \phantom{==================} 
 + V^-_1(\mu,\eta,t) + O(\sigma^{-\gamma N}), 
 $$ 
 where $$ 
 V^-_1(\mu,\eta,t)= \sum\limits_{r^2_s+q^2_s\neq 0} 
 v^-_{1,s} \mathrm{Ai}^{(m_s)}(\eta)
t^{l_s} \mu^{r_s} \ln^{q_s}\mu, 
 $$ 
 $v^-_{1,s}$ are constant coefficients. 
 
 From the explicit form the Airy function~(\ref{af}),
 it is easy to see that 
$$ 
 \int\limits_{\eta}^{+\infty} 
 \! \mathrm{Ai}(\theta)\, d\theta 
 = -\int\limits_{0}^{+\infty} 
 \frac{1}{\pi\theta} 
 \sin \left( \frac{\theta^3}{3} 
 + \eta\theta \right) d\theta 
 + \mathrm{const}. 
 $$ 
 In addition, for $\eta>0$ 
 $$ 
 \int\limits_{0}^{+\infty} 
 \frac{1}{\theta} 
 \sin \left( \frac{\theta^3}{3} 
 + {\eta\theta} \right) d\theta 
 = \int\limits_{0}^{\eta^{-1/2}} 
 + \int\limits_{\eta^{-1/2}}^{\infty} 
 = 
 $$ 
 $$ 
 =\int\limits_{0}^{\eta^{1/2}} 
 \frac{1}{z} 
 \sin \left( \frac{z^3}{3\eta^3} + z \right)
dz 
 - \int\limits_{\eta^{-1/2}}^{\infty} 
 \frac{1}{\theta (\theta^2 + \eta)} 
 d \cos \left( \frac{\theta^3}{3} 
 + {\eta\theta} \right). 
 $$ 
 The first integral on the right-hand side of this
relation as ${\eta\to +\infty}$ 
 tends to the limit value of the sine integral
$$ 
 \mathrm{Si}(+\infty)= 
 \int\limits_{0}^{+\infty} 
 \frac{\sin z}{z} dz = \frac{\pi}{2}, 
 $$ 
since $z^3/\eta^3 \leqslant \eta^{-3/2}$
on the whole integration interval.
Integrating by parts, we see that 
$$
 \int\limits_{\eta^{-1/2}}^{\infty} 
 \frac{1}{\theta (\theta^2 + \eta)} 
 d \cos \left( \frac{\theta^3}{3} 
 + {\eta\theta} \right) =O(\eta^{-1/2}). 
 $$ 
 Thus, we obtain \begin{equation}\label{iaf} 
 \int\limits_{\eta}^{+\infty} 
 \! \mathrm{Ai}(\theta)\, d\theta 
 = -\int\limits_{0}^{+\infty} 
 \frac{1}{\pi\theta} 
 \sin \left( \frac{\theta^3}{3} 
 + \eta\theta \right) d\theta 
 + \frac{1}{2}. 
 \end{equation}

 Now, let us study the integral
 $$ 
 U^-_0 (x,t) = 
 \frac{1}{\root 3 \of {3t} } 
 \int\limits_{-\sigma}^{0} 
 f(y)\, \mathrm{Ai} \left( \frac{x-y}{\root 3 \of {3t} } \right) \, dy. 
 $$ 
 Since $2t\geqslant \sigma^{\alpha/p}$ on the
set $T_{\alpha}$ 
 and  there holds the estimate $$ 
 \frac{y}{\root 3 \of {3t}} = O(\sigma^{-\gamma}), 
 \qquad 
 \gamma = \frac{\alpha}{3p}-1 > 0, 
 $$ 
for $|y| \leqslant \sigma$,
we can use Taylor's formula. 
 Then,
 $$ 
 U^-_0 (x,t)= \sum\limits_{n=0}^{N-1}
t^{-(n+1)/3} 
 \frac{(-1)^n 3^{-(n+1)/3}}{n!} 
 \mathrm{Ai}^{(n)}(\eta) 
 \int\limits_{-\sigma}^{0} 
 y^{n} f(y) dy 
 + O( \sigma^{- \gamma N}). 
 $$ 
 Let us transform the integral $\int\limits_{-\sigma}^{0}$ as follows: 
 $$ 
 \int\limits_{-\sigma}^{0} y^{n} f(y)
dy = 
 \int\limits_{-1}^{0} y^{n} f(y) dy + 
 \int\limits_{-\sigma}^{-1}\left( 
 f^-_0 y^{n} + \dots + f^-_{n} + \frac{f^-_{n+1}}{y} \right) dy + 
 $$ 
 $$ 
 + \int\limits_{-\sigma}^{-1} y^{n} \left[ 
 f(y) - f^-_0 - \dots - \frac{f^-_{n}}{y^{n}} - \frac{f^-_{n+1}}{y^{n+1}} \right]
dy = 
 $$ 
 \begin{equation}\label{is1} 
 = \int\limits_{-1}^{0} y^{n} f(y) dy + 
 \int\limits_{-\infty}^{-1}\Phi^-_n(y)
dy 
 -\int\limits_{-\infty}^{-\sigma} \Phi^-_n(y)
dy 
 - f^-_{n+1} \ln\sigma + P_n(\sigma), 
 \end{equation} 
 where $$ 
 \Phi^-_n(y) = y^{n} \left[ f(y) - \sum\limits_{m=0}^{n+1} \frac{f^-_m}{y^{m}} \right]. 
 $$ 
 From relations~(\ref{le}) we obtain 
estimates $$ 
 |\Phi^-_n(y)| \leqslant C_n y^{-2}, 
 \qquad y\leqslant -1, 
 $$ 
 $$ 
 \int\limits_{-\infty}^{-\sigma} \Phi^-_n(y)
dy = 
 \sum\limits_{m=1}^{N-1} \varphi_{n,m} \sigma^{-m} + 
 O(\sigma^{-N}), 
 \qquad \sigma\to \infty, 
 $$ 
 where $\varphi_{n,m}$ are some constants. 
 Taking into account these estimates and
substituting $\sigma=\mu{\root 3 \of {3t}}$
in (\ref{is1}), we obtain $$ 
 \int\limits_{-\sigma}^{0} y^{n} f(y)
dy = 
 I_n - \frac{f^-_{n+1}}{3} \ln t -f^-_{n+1} \ln\mu + 
 \sum\limits_{r\neq 0} a_{n,r} \mu^{r}
t^{r/3} 
 + O\left( \sigma^{-N}\right), 
 $$ 
 where $I_n$ and $a_{n,r}$ are constants. 
 Then, we have 
\begin{equation}\label{up0f} 
 U^-_0(x,t) = \sum\limits_{n=0}^{N-1} 
 t^{-(n+1)/3} \mathrm{Ai}^{(n)}(\eta) 
 [b^-_{n} + \widetilde{b}^-_{n} \ln
t] 
 +V^-_{0}(\mu,\eta,t) + O(\sigma^{-\gamma N}), 
 \end{equation} 
 $$ 
 V^-_0(\mu,\eta,t)= \sum\limits_{r^2_s+q^2_s\neq 0} 
 v^-_{0,s} \mathrm{Ai}^{(m_s)}(\eta)
t^{l_s} \mu^{r_s} \ln^{q_s}\mu, 
 $$ 
 where $b^-_{n}$, $\widetilde{b}^-_{n}$, 
 and $v^-_{0,s}$ are constant coefficients. 
 From the asymptotics of the Airy function
 \begin{equation}\label{aap} 
 \mathrm{Ai}(x) = \frac{x^{-1/4}}{2\pi} 
 \exp\left( -\frac{2}{3}x^{3/2}\right) 
 \sum\limits_{n=0}^{\infty} 
 \frac{(-1)^n \Gamma\left( \frac{3n+1}{2}\right)}{3^{2n}(2n)!} 
 x^{-3n/2}, \quad x\to +\infty, 
 \end{equation} 
 it is seen that estimates of remainders
 in formulas~(\ref{up1f}) and~(\ref{up0f}) 
 are also valid  on the set $$ 
 X^+_{\alpha} = \{ (x,t)\, : \, x>0,\ 0 <
t < |x|^{\alpha} \}, 
 $$ 
 where $$ 
 \mu = o(\eta), 
 \qquad 
 |\eta| \geqslant \frac{|x|^{1-\frac{\alpha}{3}}}{\root 3 \of 3} 
 \to \infty 
 \qquad \mbox{as} 
 \quad \sigma\to \infty. 
 $$ 

To find the dependence on $\mu$ 
of integrals $$ 
 \int\limits_{0}^{\mu} 
 \mathrm{Ai}(\eta+\theta' ) d\theta', 
 \qquad 
 \int\limits_{-\mu}^{0}\Psi_n(\theta,\eta)
d\theta, 
 \qquad 
 \int\limits_{0}^{\mu}\Psi_n(\theta,\eta)
d\theta 
 $$ 
on the set 
$$
X^-_{\alpha} = \{ (x,t)\, : \, x<0,\ 0 <t < |x|^{\alpha} \},
$$
 we use the expansion 
 \begin{equation}\label{aam} 
 \mathrm{Ai}(x) = \frac{|x|^{-1/4}}{\pi} 
 \sum\limits_{n=0}^{\infty} 
 \frac{(-1)^{[n/2]} \Gamma\left(\frac{3n+1}{2}\right)}{3^{2n}(2n)!} 
 \sin\left(\frac{2}{3}|x|^{3/2}+\frac{\pi}{4}(-1)^n \right) 
 |x|^{-3n/2} 
 \end{equation} 
as $x\to -\infty$.
Then we obtain expressions of the form 
$$ 
 \sum\limits_{r^2_s+q^2_s\neq 0} 
 b_s \mu^{r_s} (\ln\mu)^{q_s}\, 
 t^{l_s} \eta^{m_s} 
 \sin\left(\frac{2}{3}|\eta|^{3/2} \pm\frac{\pi}{4} \right), 
 $$ 
where $b_s$ are constant coefficients. 
 
 Analogously to formulas~(\ref{up1f})
and~(\ref{up0f}), we find 
\begin{equation}\label{um1f} 
 U^+_1(x,t) = 
 f^+_0\int\limits_{-\infty}^{\eta} 
 \! \mathrm{Ai}(\theta)\, d\theta 
 +\sum\limits_{n=1}^{N-1} 
 f^+_n (3t)^{-n/3} 
 \left[ \int\limits_{1}^{\infty} 
 \theta^{-n} \mathrm{Ai}(\eta-\theta)
d\theta +B^+_n(\eta) 
 \right]+ 
 \end{equation} 
 $$ 
 \phantom{==================} 
 + V^+_1(\mu,\eta,t) + O(\sigma^{-\gamma N}), 
 $$ 
 where $$ 
 V^+_1(\mu,\eta,t) = \sum\limits_{r^2_s+q^2_s\neq 0} 
 G_{1,s}(\eta) t^{l_s} \mu^{r_s} \ln^{q_s}\mu, 
 $$ 
 and smooth functions
 $$ 
 B^+_n(\eta)= \sum\limits_{r=0}^{n-2} 
 \frac{\mathrm{Ai}^{(r)}(\eta)}{r!(r-n+1)} 
 + \int\limits_{0}^{1}\Psi_n(\theta,\eta)
d\theta, 
 $$ 
 \begin{equation}\label{pinp} 
 \int\limits_{0}^{1} 
 \Psi_n(\theta,\eta) d\theta = 
 O\left( |\eta|^{-1/4+n/2} \right), 
 \qquad \eta\to \infty, 
 \end{equation} 
 are constructed similarly to functions~(\ref{bnm}), 
 \begin{equation}\label{um0f} 
 U^+_0(x,t) = \sum\limits_{n=0}^{N-1} 
 t^{-(n+1)/3} \mathrm{Ai}^{(n)}(\eta) 
 [b^+_{n} + \widetilde{b}^+_{n} \ln
t] 
 + V^+_0(\mu,\eta,t) + O(\sigma^{-\gamma N}), 
 \end{equation} 
 ${b}^+_{n}$ and $\widetilde{b}^+_{n}$ 
 are some constants,
 $$ 
 V^+_0(\mu,\eta,t) = \sum\limits_{r^2_s+q^2_s\neq 0} 
 G_{0,s}(\eta) t^{l_s} \mu^{r_s} \ln^{q_s}\mu, 
 $$ 
 $G_{0,s}(\eta)$ and $G_{1,s}(\eta)$ are smooth
functions.

 Since $\mu$ contains an arbitrary parameter~$p$ 
 and $\mu\to 0$ for $t>|x|^{h}\, (3p<h <3)$, 
 from A.R.~Danilin's lemma~\cite[Lemma 4.4]{dar}, 
 it follows that 
\begin{equation}\label{vz} 
 V^-_{0}(\mu,\eta,t)+V^-_{1}(\mu,\eta,t) 
 +V^+_{0}(\mu,\eta,t)+V^+_{1}(\mu,\eta,t)= 
 O(\sigma^{-\gamma N}). 
 \end{equation}

 Then substituting expressions~(\ref{up1f}), (\ref{up0f}), 
 (\ref{um1f}), and~(\ref{um0f}) into (\ref{ud}) 
 and using formulas~(\ref{iaf}) and~(\ref{vz}), 
 we obtain $$ 
 u(x,t)=(f^+_0-f^-_0)\int\limits_{0}^{+\infty} 
 \frac{1}{\pi\theta} 
 \sin \left( \frac{\theta^3}{3} 
 + \frac{x\theta}{\root 3 \of {3t} } \right)
d\theta 
 +\frac{f^+_0+f^-_0}{2} + 
 $$ 
 \begin{equation}\label{uxt} 
 + \sum\limits_{n=1}^{N-1} t^{-n/3} 
 \left(F_{n}(\eta) + \widetilde{F_{n}}(\eta)\ln
t \right) 
 +O(\sigma^{-\gamma N}). 
 \end{equation}

The coefficients $F_{n}(\eta)$ and $\widetilde{F_{n}}(\eta)$ 
 are linear combinations of smooth bounded functions,
derivatives of the Airy function up to the $(n-1)$th
order, inclusively, and integrals 
$$ 
 \int\limits_{-1}^{0} 
 \Psi_n(\theta,\eta) d\theta, 
 \qquad 
 \int\limits_{0}^{1} 
 \Psi_n(\theta,\eta) d\theta. 
 $$ 
Then from asymptotics~(\ref{aap}), (\ref{aam}), 
 (\ref{pinm}), and~(\ref{pinp}),
 it follows that functions $F_{n}(\eta)$
and $\widetilde{F_{n}}(\eta)$ 
 cannot grow faster than $|\eta|^{-1/4+n/2}$ 
 as $\eta\to\infty$. 
 Consequently,  the~$n$th term in expansion~(\ref{uxt}) 
 cannot grow  faster than 
 $$ 
 \ln t \, |\eta|^{-1/4} \left| \frac{x}{t}\right|^{n/2}. 
 $$ 
 Therefore, the series 
$$ 
 \sum\limits_{n=1}^{\infty} 
 t^{-n/3} \left(F_{n}(\eta) + \widetilde{F_{n}}(\eta)\ln t \right) 
 $$ 
keeps its asymptotic character for $t>|x|^{1+\delta}$ 
 with any $\delta>0$. 
 
 \vskip 3mm 
 
 Thus, we arrive at the following statement. 
 
\textbf{Theorem.} 
{\it
 If for  a locally Lebesgue integrable function~$f$ 
 the conditions 
$$ 
 f(x) = \sum\limits_{n=0}^{\infty} 
 f^{\pm}_n x^{-n}, 
 \qquad x\to \pm\infty, 
 $$ 
 are fulfilled,  then for any $\delta>0$
and $t>|x|^{1+\delta}$
 there holds the asymptotic formula 
$$ 
 \frac{1}{\root 3 \of {3t} } 
 \int\limits_{-\infty}^{+\infty} 
 f(y)\, \mathrm{Ai} \left( \frac{x-y}{\root 3 \of {3t} } \right) \, dy 
 = 
 $$ 
 $$ 
 =(f^+_0-f^-_0)\int\limits_{0}^{+\infty} 
 \frac{1}{\pi\theta} 
 \sin \left( \frac{\theta^3}{3} 
 + \frac{x\theta}{\root 3 \of {3t} } \right)
d\theta 
 +\frac{f^+_0+f^-_0}{2} 
 + \sum\limits_{n=1}^{\infty} t^{-n/3} 
 \left(F_{n}(\eta) + \widetilde{F_{n}}(\eta)\ln
t \right), 
 $$ 
as $|x| + t \to \infty$,
 where $F_{n}(\eta)$ and $\widetilde{F_{n}}(\eta)$ 
are $C^{\infty}$-smooth functions
of the self-similar variable 
$$ 
 \eta = \frac{x}{{\root 3 \of {3t}}}. 
 $$ 
}

According to formula~(\ref{uxt}),
the asymptotic series in the theorem 
 should be understood in the sense of Erdelyi
with the asymptotic sequence 
$$\{ (|x|^3 + t)^{-\gamma'n}\}_{n=1}^{\infty},$$
 where $\gamma' >0$.

\bigskip
\textbf{Acknowledgments.}
This work was supported by the Russian Foundation for Basic Research,
project no.~14-01-00322.

\end{document}